\documentclass[12pt]{article}
\usepackage{amsmath, amstext, amsgen, amsthm}
\usepackage{amsopn, amsfonts, euscript,amscd,amssymb,latexsym}
\theoremstyle{definition}

\theoremstyle{plain}
\newtheorem{Lem}{Lemma}
\theoremstyle{plain}
\newtheorem{Th}{Theorem}
\theoremstyle{plain}
\newtheorem{Cor}{Corollary}
\theoremstyle{plain}
\newtheorem{Prop}{Proposition}
\theoremstyle{definition}
\newtheorem{Rem}{Remark}
\theoremstyle{definition}
\newtheorem{Def}{Definition}
\def\id{\mathop\text{\rm id}\nolimits}
\def\Ad{\mathop\text{\rm Ad}\nolimits}
\def\Aff{\mathop\text{\rm Aff}\nolimits}
\def\Exp{\mathop\text{\rm Exp}\nolimits}
\def\Diff{\mathop\text{\rm Diff}\nolimits}
\def\Sus{\mathop\text{\rm Sus}\nolimits}
\def\Aut{\mathop\text{\rm Aut}\nolimits}
\def\Conf{\mathop\text{\rm Conf}\nolimits}
\def\GL{\mathop\text{\rm GL}\nolimits}
\def\R{\mathbb{R}}
\def\S{\mathbb{S}}

\begin{document}

\title{Attractors of Cartan
foliations}

\author{Anton S. Galaev and Nina I. Zhukova}

\maketitle

\begin{abstract}
The paper is focused on the existence problem of attractors for
foliations.  Since the existence of an attractor is a transversal
property of the foliation, it is natural to consider foliations
admitting transversal geometric structures. As transversal
structures are chosen  Cartan geometries due to their
universality. The existence problem of an attractor on a complete
Cartan foliation is reduced to a similar problem for the action of
its structure Lie group on a certain smooth manifold. In the case
of a complete Cartan  foliation with a structure subordinated to a
transformation group, the problem is reduced to the level of the
global holonomy group  of this foliation. Each countable
automorphism group preserving a Cartan geometry on a manifold and
admitting an attractor is realized as the global holonomy group of
some Cartan foliation with an attractor. Conditions on the linear
holonomy group of a leaf of a reductive Cartan foliation
sufficient for the existence of an attractor  (and a global
attractor) which is a minimal set are found. Various examples are
considered.


%foliations 53C12; 57R30;
%57Cxx; 058Axx; dg
%35B41 (Attractors)

\end{abstract}

\section{Introduction}

The study of the dynamical properties of foliations is  an actual
area. The existence of closed leaves, attractors and minimal sets
gives reach information about  the structure of a foliation. By this reason, the
problems of the existence and the structure description for
attractors and minimal sets of foliations are the central problems
in the foliation theory and  topological dynamics. There are
several nonequivalent  notions of an attractor in the theory of
dynamical systems (e.g., see \cite{Gorodetski}). Some of these
notions are equivalent  \cite{DK}. For ,,typical'' dynamical
systems  in metric sense   different notions of an attractor
coincide according to Palis's hypothesis \cite{Pal}. We use the
most general notion of an attractor for a foliation that
generalizes the notion of an attractor from \cite{SSTC}. Note that the
attractor of a foliation may be disconnected and it may contain
other attractors. This is not a case for a transitive attractor
that contains a dense leaf. Examples of transitive attractors are
attractors which are minimal sets.

In Section \ref{Sectransv} we  show that the property  of a
singular foliation to admit an attractor is transversal, i.e.,  it
is preserved under the transversal equivalence of foliations. By
this reason it is natural to investigate the influence of
different kinds of the transversal   structures of foliations on
the existence of attractors on them.  As the transversal
structures we consider Cartan geometries, since they include large
classes of geometries, e.g. Riemannian, Lorentzian (more
generally, pseudo-Riemannian), affine, conformal, projective,
transversely homogeneous, parabolic, etc.

Deroin and Kleptsyn~\cite{DK} investigated attractors of
foliations with conformal transversal structures  on compact
manifolds. The Main theorem of \cite{DK} states that for every
conformal foliation on a compact manifold either there exists a
transversely invariant measure, or there  exists  a finite number
of minimal sets equipped with probability measures, which are
attractors satisfying some properties.

The case of transversally similar foliations
considered the second author of this work in~\cite[Sec. 9]{Min}.
In \cite{ZhP}, the
  existence problem of an attractor for foliations admitting a
transversal parabolic geometry of rank one was solved. In
\cite{ZhG}, it is shown that  every non-Riemannian conformal
foliation $(M, F)$ of codimension  $q\geq 3$ admits an attractor
which is a minimal set, and the restriction of the foliation to
the basin of the attractor is a transversely conformally flat
foliation. Moreover, if the foliated manifold  $M$ is compact,
then $(M, F)$ is a $(\Conf(\S^q),\S^q)$-foliation \cite[Th.
4]{ZhL}.  Every complete non-Riemannian conformal foliation $(M,
F)$ of codimension  $q\geq 3$ admits a global attractor which is a
minimal set, and $(M, F)$ is covered by a locally trivial bundle
over the standard $q$-dimensional sphere $\S^q$ or the Euclidean
space $\R^q$ \cite[Th. 5]{ZhG}.

Note that in \cite{ZhL,ZhG} as well as in the present work we use
the methods of  local and global differential geometry, while
Deroin and Kleptsyn~\cite{DK} used  the Lyapunov exponentials and
invariant transversal  measures, including the harmonic measures.

 In Section \ref{secCartfol} we introduce the
Cartan foliations, discuss the effectivity, the completeness of
the Cartan foliations, the construction of the lifted foliation,
and the aureole foliation.

The Cartan geometry has an infinitesimal nature, consequently in
order to describe the global structure of  Cartan foliations one
should use global conditions. The most important global condition
for such foliations is the completeness. In Section
\ref{secexcret}, we use the associated singular aureole foliation
and prove a criteria (Theorem \ref{th3}) that reduces the
existence problem of an attractor for a complete Cartan foliation
of type $(G,H)$ to a similar problem for  the induced
 action of the Lie group $H$ on a certain manifold
called the basic one.

In Section \ref{S4} we study  attractors of complete Cartan
$(\Phi,N)$-foliations $(M,F)$. We prove that problems of the
existence and the structure description of attractors  (resp.
global attractors) for these foliations are equivalent to the
similar problems for the countable automorphism groups of complete
Cartan geometries on simply connected manifolds.

Next  we consider  the case of reductive Cartan foliations, i.e.,
foliations, admitting transversal reductive Cartan geometry. The
class of reductive Cartan foliations includes transversally
similar, Riemannian, Lorentzian (more generally,
pseudo-Riemannian),  reductive transver\-sal\-ly homogeneous
foliations and foliations with transversal linear connections. In
Section \ref{s3} we show that a reductive Cartan foliation admits
also a transversal linear connection; this simplifies the study of
reductive Cartan foliations.

In Section \ref{secattrminsets} we find the conditions on the
linear holonomy group of a leaf  $L$ of the foliation $(M, F)$
that are sufficient for the closure ${\mathcal M}= \overline{L}$
to be an attractor and a minimal set of this foliation. Some other
results about the linearization of the holonomy group and the
geometry around leaves of foliations obtained using various
methods Crainic, Struchiner, Weinstein, Zung, and also del~Hoyo
with Fernandes (see e.g. \cite{HF} and the references therein).
We also find the sufficient conditions
for the existence a global attractor which is a minimal set of $(M, F)$.

In Section~\ref{s7} we consider several examples.

\medskip{\noindent\bf Assumptions\,} Throughout this paper
we assume for simplicity that all manifolds and  maps are smooth
of the class $C^\infty$; in fact, the main results of the paper
are valid for foliations of the class $C^2$. All neighborhoods are
assumed to be open and all manifolds are assumed to be Hausdorff.

\medskip{\noindent\bf Notations\,} The algebra of smooth functions on a
manifold  $M$ will be denoted by  $\mathfrak{F}(M)$. Let
$\mathfrak X(N)$ denote the Lie algebra of smooth vector fields on
a manifold $N.$ If $\mathfrak M$ is a smooth distribution on $M$
and $f: K\to M$ is a submersion, then let $f^*\mathfrak M$ be the
distribution on the  manifold $K$ such that $(f^*\mathfrak M)_z=
\{X\in T_zK\,|\, f_{*z}(X)\in\mathfrak M_{f(z)}\}$, where $z\in
K$. Let $\mathfrak X_{\mathfrak M}(M)=\{X\in\mathfrak X(M)\mid
X_u\in {\mathfrak M}_u\quad\forall u\in M\}$. As usually we denote
by $P(N, H)$ the  principal $H$-bundle $P$ over the manifold~$N$.
The symbol  $\cong$ will denote the isomorphism of objects in the
corresponding category.

\section{Attractors of foliations and transversality}\label{Sectransv}

In this section we give a definition of an attractor of a singular
foliation in the sense of Stefan and Sussmann and show that the
property of a foliation to admit an attractor is transversal. Most
of the results of this paper are obtained for smooth foliations.
We use singular foliations in the proof of Theorem~\ref{th3}

\begin{Def}\label{d1} Let $(M, F)$ be a singular foliation. A subset of a manifold
 $M$ is called saturated if it is a union of leaves of this foliation.
A nonempty closed saturated subset $\mathcal M$ of $M$ is called
{\it an attractor} of  $(M, F)$ if there exists an open saturated
neighbourhood $\mathcal U = \mathcal U(\mathcal M)$ of the set
$\mathcal M$ such that the closure of every leaf from $\mathcal
U\setminus{\mathcal M}$ contains the set $\mathcal M$, i.e., if
$\overline{L}\supset\mathcal M$ $\forall L\subset\mathcal
U\setminus{\mathcal M}.$ The neighbourhood $\mathcal U$ is
uniquely determined by this condition and it is called {\it the
basin of this attractor}; we denote it by ${\cal B}(\mathcal M)$.
If in addition ${\cal B}({\mathcal M}) = M$, then the attractor
$\mathcal M$ is called {\it global}.
\end{Def}

\begin{Def} Two smooth singular foliations $(M_1,F_1)$ and $(M_2,F_2)$ are called
transversally equivalent if there
exists a smooth singular foliation $(\mathbb M, \mathbb F)$ and two
submersions $p_1: \mathbb M\to M_1$ and $p_2: \mathbb M\to M_2$
such that
$$\mathbb F = \{p_1^{-1}(L_\alpha) |\, L_\alpha\in F_1\} = \{p_2^{-1}(L_\beta) |\, L_\beta\in F_2\}.$$
\end{Def}
This notion generalizes the notion of the transversal equivalence
for smooth foliation in the sense of Molino \cite[Def. 2.1]{Mo}.
It can be checked directly that the transversal equivalence of
singular foliations is an equivalence relation.

A property of singular foliations is said to be {\it transversal} if it
is preserved under the transversal equivalence.
%Recall that a {\it saturated set} of a foliated manifold is a union of leaves.

\begin{Prop}\label{p1} The existence of an attractor is
a transversal property of a foliation.
\end{Prop}

\vspace{3mm} {\bf Proof of Proposition~\ref{p1}.} Since
transversally equivalent singular foliations have the common leaf
space, it is sufficient to characterize  the existence of an
attractor in terms of the topology of the leaf space $M/F$ of a
singular foliation $(M, F)$.

It is well known  that the projection $f: M\to M/F$ is an open map
and the closure in $M$ of any saturated subset is again a
saturated subset. Due to this we observe that a singular foliation
$(M, F)$ has an attractor if and only if there exists a nonempty
closed subset $\widetilde{\mathcal M}\subset M/F$ and its open
neighbourhood $\widetilde{\mathcal U}$ in $M/F$ such that the
closure  $\overline{\{z\}}$ of any one-point set
$\{z\}\subset\widetilde{\mathcal U}\setminus\widetilde{\mathcal
M}$ satisfies the inclusion $\widetilde{{\mathcal M}}\subset
\overline{\{z\}}.$ Note that the set ${\mathcal M} =
f^{-1}({\widetilde{\mathcal M}})$ is an attractor of the singular
foliation $(M, F)$. This proves the proposition. \qed

\begin{Def}\label{d2} Let  $H$ be a Lie group with a  smooth action on a manifold $W$.
A nonempty closed union $\mathcal K$ of orbits of $H$ is called an
attractor of the action if there exists an open  invariant
neighbourhood $\mathcal V = \mathcal V(\mathcal K)$ of the set
$\mathcal K$ such that the closure of every orbit from $\mathcal
V\setminus{\mathcal K}$ contains the set $\mathcal K$. If
${\mathcal V} = W$, then the attractor $\mathcal K$ is called {\it
global}.
\end{Def}

The following proposition can be proved in the same way as
Proposition~\ref{p1} (for connected Lie groups this proposition
follows from Proposition~\ref{p1}).

\begin{Prop}\label{p2} The existence of an attractor of a smooth action of
a Lie group on a manifold is a transversal property.
\end{Prop}

\section{Cartan foliations and  the associated
constructions}\label{secCartfol}
\subsection{Cartan geometries}
Let us first recall  the definition of a Cartan geometry
\cite{C-S,CrS}. Let $G$ be a Lie group and  $H$ be a closed
subgroup of $G$. Denote by $\mathfrak{g}$ and $\mathfrak{h}$ the
Lie algebras of the Lie groups $G$ and $H$, respectively. Let $N$
be a smooth (not necessary connected) manifold. {\it A Cartan
geometry} on $N$ of type $(G,H)$ (or $\mathfrak g/\mathfrak h$) is
a principal right $H$-bundle $P(N,H)$ with the projection
$p:P~\rightarrow~N$ and together with a $\mathfrak{g}$-valued
$1$-form $\beta$ on $P$ satisfying the following  conditions:
\begin{enumerate} \itemsep=0pt
\item[({$c_1$})] the map $\beta_{w}:T_{w}P\rightarrow \mathfrak{g}$ is an isomorphism of
the vector spaces for every $w\in P$;
\item[($c_{2}$)] $R^{*}_{h}\beta=\Ad_{G}(h^{-1})\beta$ for all $h\in
H$, where $\Ad_{G}:G\rightarrow \GL(\mathfrak{g})$  is the adjoint
representation of the Lie group  $G$ on its Lie algebra
$\mathfrak{g}$;
\item[($c_{3}$)]  $\beta({A^{*}})=A$ for any $A\in\mathfrak{h}$, where $A^{*}$
is the fundamental vector field defined by the element $A$.
\end{enumerate}
 This Cartan geometry is denoted by $\xi=(P(N,H),\beta)$.
The pair $(N,\xi)$ is called a {\it Cartan manifold}.

\begin{Def}\label{d6} Let $G/H$ be a homogeneous space and let $G$ act on $G/H$
by  left  translations. Denote by $\mathfrak{g}$ and
$\mathfrak{h}$ the Lie algebras of the Lie groups $G$ and $H$,
respectively. If there exists an $\Ad_G(H)$-invariant vector
subspace $V$ of $\mathfrak g$ such that
$$\mathfrak g = \mathfrak h \oplus V,$$
then the homogeneous space $G/H$ is called reductive. A Cartan geometry $\xi=(P(N,H),\beta)$
of type $(G,H)$, where $G/H$ is a reductive homogeneous space, is called {\it a
reductive Cartan geometry}.
\end{Def}

Cartan manifolds form a category, where morphisms of two Cartan
manifolds $\xi=(P(N,H),\beta)$ and $\xi'=(P'(N',H),\beta')$ of the
same type $(G,H)$ are morphisms of the principle bundles
$\Gamma:P(N,H)\to P'(N',H)$ satisfying the  condition
$\Gamma^*\beta'=\beta$.

\subsection{Cartan foliations}\label{ss3.2}
A foliation $(M,F)$ is said to be {\it a Cartan foliation} of type
$(G,H)$ if it admits a Cartan geometry of type $(G,H)$ (or
$\mathfrak g/\mathfrak h$) as a transversal  structure. More
precisely, this means the following. Let $\xi=(P(N,H),\beta)$ be a
Cartan geometry. A Cartan foliation $(M,F)$ may be defined using
an $(N,\xi)$-\textit{cocycle}, i.e. a family $\{U_i,
f_i,\Gamma_{ij}\}_{i,j\in J}$ satisfying  the following
properties:
\begin{enumerate}
\item [(i)] $\{U_i\}_{i\in J}$ is an open covering of the manifold  $M$ by connected subsets
$U_i$, and  $f_i\colon U_i \to V_i \subset N$ are submersions with
connected leaves;
\item[(ii)] if $U_i \cap U_j\neq\emptyset$, $i,j\in J$, then there exists an isomorphism  $\Gamma_{ij}$
of the Cartan geometries $\xi_{f_j(U_i \cap U_j)}$ and
$\xi_{f_i(U_i \cap U_j)}$ induced on $f_j(U_i \cap U_j)$ and
$f_i(U_i \cap U_j)$, respectively, such that its projection
$\gamma_{ij}$ satisfies the equality $f_i=\gamma_{ij} \circ f_j$;
\item  [(iii)] $\Gamma_{ij}\circ \Gamma_{jk} = \Gamma_{ik}$ for $U_i\cap U_j \cap
U_k\neq\emptyset$; moreover, $\Gamma_{ii}=\id_{P_{f_i(U_i)}}$.
\item  [(iv)] the cocycle $\{U_i,
f_i,\gamma_{ij}\}_{i,j\in J}$ defines the foliation $(M,F)$.
\end{enumerate}

One says also that the Cartan foliation $(M,F)$ satisfying the
above properties is modelled on the Cartan geometry $\xi$ of type
$(G,H)$.

\subsection{The lifted foliation}\label{ss3.2A}

We will use the construction of the lifted foliation $(\cal R,
\cal F)$ for a Cartan foliation $(M,F)$ from \cite{Min}; it
generalizes a similar construction for a Riemannian foliation from
\cite{Mo}. For a given Cartan foliation $(M,F)$ of type $(G,H)$
one may construct a principle $H$-bundle ${\cal R}(M,H)$  (called
a foliated bundle) with a projection $\pi:{\cal R}\to M$, an
$H$-invariant transversely parallelizable foliation $({\cal
R},{\cal F})$ such that  $\pi$ is a morphism of $({\cal R},{\cal
F})$ into $(M,F)$ in the category of foliations; moreover, there
exists a $\mathfrak g$-valued 1-form $\omega$ on ${\cal R}$ having
the following properties:

(i) $\omega(A^*)=A$ for any $A\in \mathfrak h,$ where $A^*$ is the
fundamental vector field corresponding to $A$;

(ii) $R_a^*\omega=\Ad_G(a^{-1})\omega$ $\forall a\in H$;

(iii) for any $u\in {\cal R}$, the map $\omega_u: T_u{\cal R}\to
\mathfrak g$ is surjective with the kernel $\ker \omega = T{\cal
F}$, where $T{\cal F}$ is the tangent distribution to the
foliation $(\cal R, \cal F)$;

(iv) the Lie derivative $L_X\omega$ is zero for any vector field $X$ tangent to the leaves of $(\cal R, \cal F).$

The foliation $(\cal R, \cal F)$ is called the {\it lifted
foliation}. The restriction $\pi_{\mathcal L}: {\mathcal L}\to L$
of $\pi$ to a leaf $\mathcal L$ of $(\cal R, \cal F)$  is a
covering map onto the corresponding leaf $L$ of $(M, F)$. If $\cal
R$ is disconnected, then we consider a connected component of
$\cal R$.

\subsection{Effectivity of transversal Cartan geometries}\label{ss2.3}
Let us recall several standard definitions.

\begin{Def} A pair of Lie algebras $(\mathfrak g,\mathfrak h)$, where
$\mathfrak h$ is subalgebra of $\mathfrak g$, is called effective
if the maximal ideal of $\mathfrak g$ belonging to $\mathfrak h$
is zero.
\end{Def}

\begin{Def} A Cartan geometry of the type $(G,H)$ is said to be
{\it effective} if the group $G$ acts effectively on  $G/H$, or in
other words, if the maximal normal subgroup of $G$ belonging to
$H$ is trivial.
\end{Def}

Note that the effectivity of a pair of Lie groups $(G,H)$, where
$H$ is a closed subgroup of $G$, implies the effectively of the
pair $(\mathfrak g,\mathfrak h)$ of their Lie algebras. It is
known \cite[Prop. 1]{Min} that if $(M, F)$ admits an  ineffective
Cartan geometry, then $(M, F)$ admits also an effective
 Cartan geometry. Therefore without loss of generality we assume
further in this work that all Cartan foliations are modelled on
effective Cartan geometries  if the contrary is not indicated.

In the case of an effective  Cartan geometry, the definition of a
Cartan foliation  from the previous section is equivalent to the
definition of a Cartan foliation by Blumenthal \cite{Bl}.

\subsection{Completeness of Cartan foliations}\label{seccomplit}
Let $(M, F)$ be a foliation. A $q$-dimensional smooth distribution
$\mathfrak M$ on the manifold $M$ is called transversal   to the
foliation $(M, F)$ if the equality $T_{x}M= T_x{F}\oplus{\mathfrak
M}_x$ holds for all $x\in M$. One may identify a transversal
distribution $\mathfrak{M}$ with the vector quotient bundle
$TM/TF$.

Let $(M, F)$ be a Cartan foliation of codimension $q$ and let
$\mathfrak M$ be a transversal  $q$-dimensional distribution.
Denote by $\widetilde{\mathfrak M}$ the induced distribution
$\pi^*\mathfrak M$ on $\cal R.$ Let $\omega$ be the $\mathfrak
g$-valued $1$-form on $\cal R$ defined in the previous section.
The  Cartan foliation $(M, F)$ is said to be $\mathfrak M$-complete
if every vector field $X\in\frak X_{\widetilde{\mathfrak M}}(\cal
R)$ satisfying the condition $\omega(X)=c$, $c = const\in\mathfrak
g,$ is complete. A  Cartan foliation $(M, F)$ is called {\it
complete} if there exists a transversal   distribution $\mathfrak
M$ such that $(M, F)$ is $\mathfrak M$-complete.

\subsection{The associated aureole foliation}\label{secaureol}
Consider a complete Cartan foliation $(M, F)$ and its lifted
foliation $(\mathcal R,\mathcal F)$. Let $\pi: {\mathcal R}\to M$
be the projection of the corresponding $H$-bundle. Since
$(\mathcal R,\mathcal F)$ is a complete transversally
parallelizable foliation, the closures $\overline{\mathcal L}$ of
its leaves $\cal L$ form a simple foliation $({\mathcal
R},\overline{\mathcal F})$ which leaves are fibres of a locally
trivial bundle $\pi_B: {\mathcal R}\to W$ \cite[Th.4.2]{Mo}.  The
manifold $W$ is called the {\it basic manifold associated to the
foliation}  $(M,F)$. The image ${\cal O}(L)=
\pi(\overline{\mathcal L})$, where ${\mathcal L}\in {\mathcal F}$,
is called the {\it aureole}  of the leaf $L = \pi({\mathcal L})$
of $(M, F)$. The aureole ${\mathcal O}(L)$ of a leaf $L = L(x)$ is
also denoted ${\cal O}(x)$.

\vspace{3mm} {\bf Theorem}  \cite[Th. 2]{Min}. {\it The set of all
aureoles ${\mathcal O}= \{\pi(\overline{\mathcal L})\,|\,
{\mathcal L}\in{\mathcal R}\}$ of a complete Cartan foliation
$(M,F)$ is a smooth singular foliation $(M, {\mathcal O})$ that
has the following properties:

1) the leaf $L(y)$ of $(M,F)$ is dense in ${\mathcal O}(x)$, $x\in M$,
for every point $y\in {\mathcal O}(x)$;

2) $L(x)\subset{\mathcal O}(x)\subset\overline{L(x)}$ $\forall
x\in M$, where $\overline{L(x)}$ is the closure of $L(x)$ in $M$.}
%\end{Theorem}

\vspace{3mm} The foliation $(M,{\mathcal O})$ defined in the above
theorem is called {\it the aureole foliation} associated with $(M,
F)$. The map
\begin{equation}\label{eq1B}
\Phi^W: W\times H\to W: (w,a)\mapsto \pi_b(R_a(u))\,\,\,\, \forall\, (w,a)\in W\times H,
u\in\pi_b^{-1}(w), %\nonumber
\end{equation}
defines an action of the Lie group $H$ on the basic manifold $W =
{\mathcal R}/\overline{\mathcal F}$, and the orbit space $W/H$ is
homeomorphic to the leaf space $M/\cal O$ of the aureole foliation
$(M,\mathcal O)$, i.e. $W/H\cong M/\cal O$.

\section{A criteria for the existence  of an
attractor for a complete Cartan foliation}\label{secexcret}

\begin{Th}\label{th3}
Let $(M,F)$ be a complete Cartan foliation of type $(G,H)$. The
foliation  $(M,F)$ admits an attractor  (or a global attractor) if
and only if the induced action of the Lie group $H$ on the basic
manifold has an attractor  (or a global attractor).
\end{Th}

{\bf Proof of Theorem~\ref{th3}.} Assume that a complete Cartan
foliation $(M, F)$ admits an attractor $\mathcal M$ with the basin
${\mathcal B} = {\mathcal B}(\mathcal M).$ Let us show that
$\mathcal M$ is an attractor for the associated singular aureole
foliation $(M,\mathcal O)$ with the same basin $\mathcal B$. Pick
a point $x\in{\mathcal B}\setminus\mathcal M$. Let $L = L(x)$,
then $L\cap{\mathcal M} = \emptyset$. By the definition of an
attractor, we get
%\begin{equation}\label{q3}
$${\mathcal M}\subset\overline{L}.$$
%\end{equation}
Consider the aureole ${\mathcal O}(L)$. According to the above
theorem, $\overline{L} = \overline{{\mathcal O}(L)}$, hence
\begin{equation}\label{eq2}
{\mathcal M}\subset\overline{{\mathcal O}(L)}.
\end{equation}
First we check that
\begin{equation}\label{eq3}
{\mathcal M}\cap{\mathcal O}(L) = \emptyset.
\end{equation}
Suppose that ${\mathcal M}\cap{\mathcal O}(L)\neq\emptyset.$ Since
both ${\mathcal M}$ and ${\mathcal O}(L)$ are saturated sets,
there exists a leaf $L_\alpha\subset{\mathcal M}\cap{\mathcal
O}(L)$. Therefore we get the following chain of relations
$$\overline{{\mathcal O}(L)} =
\overline{L_\alpha}\subset\overline{{\mathcal M}\cap{\mathcal
O}(L)}\subset {\mathcal M}\cap\overline{{\mathcal
O}(L)}\subset{\mathcal M},$$ which together with (\ref{eq2})
implies $\overline{{\mathcal O}(L)} = \mathcal{M}$. Consequently,
$L\subset\overline{L} = \overline{{\mathcal O}(L)}\subset\mathcal
M$. This contradicts  $L\cap{\mathcal M} = \emptyset$. Thus
(\ref{eq3}) holds true.

Pick a leaf $L_\beta\subset{\mathcal O}(L)$, then
$$\overline{L_\beta} = \overline{{\mathcal O}(L)} =
\overline{L}\supset\mathcal{M}. $$ Therefore,
$\overline{L_\beta}\supset{\mathcal M}$. Consequently,
$L_\beta\cap{\mathcal B}\neq\emptyset$. Since ${\mathcal B}$ is a
saturated set, it holds  $L_\beta\subset{\mathcal B}$. This and
(\ref{eq3}) imply that ${\mathcal O}(L)\subset{\mathcal
B}\setminus{\mathcal M}$. Consequently $\mathcal M$ is an
attractor of the aureole foliation $(M,\mathcal O)$ with the same
basin ${\mathcal B} = {\mathcal B}({\mathcal M})$.

From  Propositions~\ref{p1} and \ref{p2}  it follows that the
induced action of the Lie group $H$ on the basic manifold $W$ has
the attractor ${\mathcal K} = \pi_B(\pi^{-1}(\mathcal M))$, and
${\mathcal B}({\mathcal K}) = \pi_B(\pi^{-1}(\mathcal B(\mathcal
M)))$ is its basin.

Conversely, assume that the induced action of the Lie group $H$ on
the basic manifold $W$ has an attractor ${\mathcal K}$. Since
$M/{\mathcal O}\cong W/H$, by  Propositions~\ref{p1} and \ref{p2},
the aureole foliation $(M,\mathcal O)$ admits the attractor
$\mathcal M= \pi(\pi_B^{-1}({\mathcal K}))$. Let $\mathcal B =
\mathcal B(\mathcal M)$ be its basin. Consider any leaf
$L_\beta\subset{\mathcal B}\setminus\mathcal M$. According to the
above theorem, $\overline{L_\beta} = \overline{{\mathcal
O}(L_\beta)}\supset\mathcal M$. Therefore $\mathcal M$ is an
attractor of $(M, F)$ with the same basin $\mathcal B(\mathcal
M)$. Since ${\mathcal B}({\mathcal M}) = M$ if and only if
${\mathcal B}({\mathcal K}) = W$, $\mathcal M$ is a global
attractor of the foliation $(M, F)$ if and only if $\mathcal K$ is
a global attractor of the group $\Psi$. \qed

\section{Attractors of $(\Phi,N)$-foliations}\label{S4}

\subsection{$(\Phi,N)$-manifolds and $(\Phi,N)$-foliations}
Let $N$ be a connected manifold and let $\Phi$ be a group of
diffeomorphisms   of $N$. One says that a group $\Phi$ acts
quasi-analytically  on $N$ if for any open subset $U$ of $N$ the
condition $\phi|_{U} = \id_U$ implies $\phi =  \id_N$. In this
section we assume that a group $\Phi$ of diffeomorphisms of a
manifold $N$ acts quasi-analytically.

\begin{Def} A foliation $(M, F)$ defined by an $N$-cocycle $\{U_i,f_i,\gamma_{ij}\}_{i,j\in J}$
is called a $(\Phi,N)$-foliation if for every $U_i\cap
U_j\neq\emptyset$, $i,j\in J$, there exists an element
$\phi\in\Phi$ satisfying the equality $\phi|_{f_j(U_i\cap U_j)} =
\gamma_{ij}$.
\end{Def}

\begin{Def} A manifold $B$ is called a $(\Phi,N)$-manifold if its natural zero-dimensional
foliation $(B, F^0)$ is a $(\Phi,N)$-foliation.
\end{Def}
We emphasize that a group of automorphisms of a Cartan manifold $(N,\xi)$ acts quasi-analytically on $N$.

\begin{Def} If $\Phi$ is a subgroup of the Lie group of all automorphisms of a Cartan
manifold $(N,\xi)$, then a $(\Phi,N)$-foliation is called a
Cartan $(\Phi,N)$-foliation.
\end{Def}

\begin{Th}\label{th4} Let $(M, F)$ be a complete Cartan $(\Phi,N)$-foliation. Then
\begin{itemize} \item[(i)] there exists a regular covering map $\kappa: \widehat M \to
M$ such that the induced  foliation $(\widehat M, \widehat F)$, $
\widehat F = \kappa^*F$, consists of the fibres of a locally
trivial bundle  $r:\widehat M\to B$, where $(B,\xi)$ is a simply
connected Cartan $(\Phi,N)$-manifold with a complete Cartan
geometry $\xi$;
\item[(ii)] an epimorphism  $\chi: \pi_1(M)\to\Psi$ of the fundamental group $\pi_1(M)$ onto a subgroup
$\Psi$ of the automorphism group $\Aut(B,\xi)$ of the Cartan
manifold $(B,\xi)$ is defined in such a way that  $\Psi$ is
isomorphic to the deck  transformation group of the covering
$\kappa: \widehat M \to M$;
\item[(iii)] for all points $y\in M$ and $z\in\kappa^{-1}(y)$, the restriction
$\kappa|_{\widehat{L}}: \widehat{L}\to L$ to the leaf $\widehat{L}
= \widehat{L}(z)$ of the foliation $(\widehat{M},\widehat{F})$ is
a regular covering map onto the leaf $L = L(y)$ of the foliation
$(M, F)$,  the group of the deck transformations is isomorphic to
the stationary sub\-group $\Psi_v$ of $\Psi$ at the point $v=
r(z)\in B$, and $\Psi_v$ is isomorphic to the holonomy group
$\Gamma(L,y)$ of the leaf $L$.
\end{itemize}

Moreover, $(M, F)$ has an attractor (resp., a global attractor)
$\mathcal M$ if and only if the group  $\Psi$ has an attractor
(resp., a global attractor) $\mathcal K$, and ${\mathcal M} =
\kappa(r^{-1}(\mathcal K))$. Besides, $\mathcal M$ is a minimal
set of the foliation $(M, F)$ if and only if $\mathcal K$ is a
minimal set of the group $\Psi$.
\end{Th}

\begin{Def}\label{d3A} The group $\Psi$ appearing in Theorem~\ref{th4} is called the global
holonomy group of the foliation $(M, F)$.
\end{Def}

\begin{Cor}\label{c1} The transversal structure of a global attractor of a foliation
$(M, F)$ satisfying Theorem~\ref{th4} is completely determined by
the structure of the corresponding global attractor of its global
holonomy group $\Psi$.
\end{Cor}

\begin{Th} \label{th5} Let  $(B,\xi)$ be a simply connected  Cartan manifold.
Let $\Psi$ be a countable subgroup of the Lie group $\Aut(B,\xi)$
of all automorphisms of $(B,\xi)$. Suppose that $\Psi$ has an
attractor (resp., a global attractor). Then $\Psi$  may be
realized as the global holonomy group of a certain Cartan
$(\Phi,N)$-foliation admitting an attractor (resp., a global
attractor).
\end{Th}

Theorems~\ref{th4} and~\ref{th5} show that the problems of the
existence and the structure description of attractors (resp.
global attractors) of complete Cartan $(\Phi,N)$-foliations  are
equivalent to the similar problems for countable automorphism
groups of complete Cartan geometries on simply connected
manifolds.

\subsection{Ehresmann connection for foliations}
The notion of an Ehresmann connection for  foliations was
introduced by Blumenthal and Hebda in~\cite{BH}. We use the
terminology from \cite{Min}. Let $(M, F)$ be a smooth
 foliation of codimension $q\geq 1$ and  $\mathfrak{M}$ be a
 $q$-dimensional  transversal distribution on $M$. All maps considered here are assumed to be piecewise smooth.
The curves in the leaves of the foliation are called vertical;
the distribution $\mathfrak{M}$ and its  integral curves are
called horizontal.

A map $H:I_1\times I_2\rightarrow M,$ where
 $I_1=I_2=[0,1]$, is called {\it a vertical-horizontal homotopy} if
for each fixed  $t\in I_2$, the curve $H_{|I_1\times \{t\}}$ is
horizontal, and for each fixed  $s\in I_1$, the curve
$H_{|\{s\}\times I_2}$ is vertical. The pair of curves
$(H_{|I_1\times \{0\}},H_{|\{0\}\times I_2})$ is called the base
of $H$.

A pair of curves  $(\sigma,h)$ with a common starting point
$\sigma(0) = h(0)$, where $\sigma:I_1\rightarrow M$
 is a horizontal curve, and $h:I_2\rightarrow M$ is a vertical curve, is called admissible.
If for  each admissible pair of curves $(\sigma,h)$ there exists a
vertical-horizontal homotopy with the base $(\sigma,h)$, then the
distribution  $\mathfrak M$ is called an Ehresmann connection for
the  foliation $(M, F)$. Note that there exists at most one
vertical-horizontal homotopy with a given base. Let $H$ be a
vertical-horizontal homotopy with the base  $(\sigma,h)$. We say
that $\widetilde \sigma= H|_{I_1\times \{1\}}$ is the result of
the  translation of the horizontal curve $\sigma$ along the
vertical curve $h$ with respect to the Ehresmann connection
$\mathfrak M$. Similarly the curve $\widetilde h=H_{|\{1\}\times
I_2}$ is called the translation of the curve $h$ along $\sigma$
with respect to $\mathfrak M$. We use the denotation
$\sigma\stackrel {h}{\mapsto }\widetilde{\sigma}$ and $h\stackrel
{\sigma}{\mapsto}\widetilde{h}$.

\subsection{Proof of Theorem~\ref{th4}}
Let $(M, F)$ be a complete Cartan $(\Phi,N)$-foliation of
codimension $q$. Then there exists a transversal $q$-dimensional
distribution $\mathfrak M$ such that $(M, F)$ is $\mathfrak
M$-complete. According to \cite[Prop. 2]{Min},  $\mathfrak M$ is
an Ehresmann connection for $(M, F)$. Applying \cite[Th. 2]{ZhG}
to the $(\Phi,N)$-foliation $(M, F)$, we see that there exists a
regular covering $\kappa: \widehat{M}\to M$ such that the induced
foliation $(\widehat{M},\widehat{F})$, $\widehat{F}=\kappa^*F$, is
made up of fibres of the locally trivial bundle $r: \widehat{M}\to
B$ over a simply connected smooth manifold $B$. Besides, there is
the induced group $\Psi$ of diffeomorphisms of $B$ and an
epimorphism
$$\chi: \pi_1(M,x)\to\Psi$$
of the fundamental group $\pi_1(M,x)$ of $M$, $x\in M$, onto
$\Psi$. Further, the group of deck transformations of the covering
$\widehat{M}$ is isomorphic to the group $\Psi$. Note that the
foliation $(\widehat{M},\widehat{F})$ is a Cartan foliation, and
it is $\widehat{\mathfrak M}$-complete with respect to the induced
distribution $\widehat{\mathfrak M}$, where $\widehat{\mathfrak M}
= \kappa^*\mathfrak M.$ Observe that the transversal Cartan
geometry of the foliation $(\widehat{M},\widehat{F})$ induces a
complete Cartan geometry $\eta$ on $B$. The group $\Psi$ is a
countable subgroup of $\Aut(B,\eta)$ of  the Lie group of all
automorphisms of $(B,\eta)$. We have proven the statements $(i)$
and $(ii)$ of Theorem~\ref{th4}. The statement $(iii)$ of
Theorem~\ref{th4} follows from the similar statement
\cite[Th.~2]{ZhG}.

Assume that there exists an attractor (resp., a global attractor)
$\mathcal M$ of the foliation $(M, F)$. It is easy to check  that
$\mathcal K= r(\kappa^{-1}(\mathcal M))$ is an attractor (resp., a
global attractor)  of the group~$\Psi$. Conversely, let $\mathcal
K$ be an attractor (resp., a global attractor) of the group
$\Psi$. It is easy to see that $\mathcal M= \kappa(r^{-1}(\mathcal
K))$ is an attractor (resp., a global attractor) of the foliation
$(M, F)$. Finally,  $\mathcal K$ is a minimal set of the group
$\Psi$ if and only if $\mathcal M = \kappa(r^{-1}(\mathcal K))$ is
a minimal set of the foliation $(M, F)$. \qed

\subsection{The suspended  foliation}
Let  $B$ and $T$ be connected smooth manifolds of dimensions
 $n-q$ and $q$, respectively. Let $\rho:
\pi_1(B,b)\to \Diff(T)$ be a homomorphism  from the fundamental
group $G = \pi_1(B,b)$ to the  group of diffeomorphisms  of the
manifold $T$. We consider the universal covering space $\widehat
B$ of $B$ as a right $G$-space. Let us define the left  action of
the group $G$ on the product $\widehat{B}\times T$ by the rule
$$\Phi: G\times\widehat{B}\times T\to\widehat{B}\times T,\quad (g,
(\widehat{b},t))\mapsto (\widehat{b}.g^{-1}, \rho(g) t),$$ where
$(\widehat{b},t)\in\widehat{B}\times T$. Then we obtain a smooth
 $n$-dimension quotient manifold $M =
\widehat{B}\times_G T$ with a foliation $F$ of codimension $q$.
The leaves of the foliation $(M, F)$ are images of the leaves of
the foliation $F_{tr} = \{\widehat{B}\times\{t\} |\, t\in T\}$
under the  quotient map $f: \widehat{B}\times T\to M$, which is a
regular covering. The foliation $(M, F)$ is called  {\it the
suspension} and it is denoted by  $\Sus(T,B,\rho)$. One says that
$(M, F)$ is obtained from the  {\it suspension of the
homomorphism} $\rho$.

The images of the leaves of the  foliation $F^t_{tr} =
\{\{\widehat{b}\}\times T |\, \widehat{b}\in\widehat{B}\}$ on the
product manifold $\widehat{B}\times T$ form a locally trivial
bundle $p: M\to B$, which is transversal to the foliation $(M,
F)$.

\subsection{Proof of Theorem~\ref{th5}}
Let $(T,\eta)$ be a simply connected $q$-dimensional Cartan
manifold and let $\Aut(T,\eta)$ be the Lie group of all
automorphisms of $(T,\eta)$. Assume that  $\Psi$ is a countable
subgroup of the  group $\Aut(T,\eta)$ and assume that $\Psi$
admits an attractor ${\mathcal K}\subset T.$

First we suppose that $\Psi$ has a finite set of generators
$\{\psi_1,...,\psi_m\}$. Denote by $\S^2_m$ the $2$-dimensional
sphere with $m$ handles. As it is known, the fundamental group of
$\S^2_m$ may be represented in the form
$$\pi_1(\S^2_m) =
<a_1,\dots, a_m,b_1,\dots, b_m| a_1b_1a_1^{-1}b_1^{-1}...
a_mb_ma_m^{-1}b_m^{-1}=1>.$$ Let $B=\S^2_m$ and define the
homomorphism $\rho : \pi_1 (B, b)\to \Aut(T,\eta)$ by the
conditions
$$\rho(a_i)=\psi_i,\quad \rho(b_i)=\id_{T}, \quad i=1,\dots,m,$$ here
$\id_{T}$ is the neutral element of the group $\Psi$. Then we
consider the suspended foliation $(M, F)=\Sus(T, B, \rho)$. The
foliation $(M, F)$ is a Cartan foliation of codimension $q$
covered by the locally trivial bundle $\R^2\times T\rightarrow T$,
and $\Psi$ is its global holonomy group. The manifold $M$ is the
total space of the locally trivial bundle  $p: M\to B$ with the
standard leaf $T$ over the base $B$. From the compactness of the
manifolds $T$ and $B$ follows the compactness of $M$.

Suppose now that the group $\Psi\subset \Aut(T,\eta)$ has a
countable set of generators  $\{\psi_i\,|\, i\in\mathbb N\}$. Let
$T_\infty$ be the plain with the pitched  countable subset
$\{(n,0)\,|\, n\in\mathbb{N}\}$. Then, $$\pi_1(T_\infty)=<a_n
\,|\,n\in\mathbb{N}>.$$ The assignment
$$\rho_\infty(a_n)=\psi_n,\quad n\in\mathbb{N},$$ defines the
homomorphism $$\rho_{\infty}: \pi_1 (B, b)\to \Aut(T,\eta).$$ The
suspended foliation $(M, F)$=$\Sus(T,B,\rho_\infty)$  is a Cartan
foliation of codimension $q$ with the global holonomy group
$\Psi$.

By the assumption, the group $\Psi$ has an attractor ${\mathcal
K}\subset T.$ Let us consider the regular covering $f: \R^2\times
T\to M$ and the projection onto the second factor $r: \R^2\times
T\to T$. Then ${\mathcal M}= f(r^{-1}(\mathcal K))$ is an
attractor of the foliation $(M, F)$. It is easy to see that
${\mathcal M}$ is a global attractor (resp., a minimal set) of
$(M, F)$ if and only if ${\mathcal K}$ is a global attractor
(resp., a minimal set) of $\Psi$. \qed

\section{Reductive Cartan foliations as foliations \\with  transversal linear connections}
\label{s3}
\subsection{Foliations with transversely linear connection}
Let $(N^{(i)},\nabla^{(i)})$, $i = 1,2$, be manifolds with linear
connections $\nabla^{(i)}$. A diffeomorphism  $f:N^{(1)}\to N^{(2)}$ is called an
isomorphism of the connections $\nabla^{(1)}$ and $\nabla^{(2)}$
if
$$f_*(\nabla^{(1)}_XY)=\nabla^{(2)}_{f_*X} f_*Y$$ for all vector fields
$X,Y\in\mathfrak X(M^{(1)})$, where $f_*$ is the differential of  $f$.

\begin{Def}~\label{d12} Suppose that an $N$-cocycle
$\left\{U_{i},f_{i},\gamma_{ij}\right\}_{i,j\in J}$ defines the
foliation $(M, F)$. If on the manifold  $N$ a linear connection
$\nabla$ is given such that each local diffeomorphism
$\gamma_{ij}$ is an isomorphism of the linear connections induced
by $\nabla$ on open subsets $f_{i}(U_{i}\cap U_{j})$ and
$f_{j}(U_{i}\cap U_{j}),$ then $(M, F)$ is called {\it a foliation
with a transversely linear connection} given by the
$(N,\nabla)$-cocycle
$\left\{U_{i},f_{i},\gamma_{ij}\right\}_{i,j\in J}$. It is said
that  $(M, F)$ is modeled on the manifold with the linear
connection $(N,\nabla)$. We stress that the connection $\nabla$ on
$N$ may have a nonzero torsion.
\end{Def}

\begin{Rem}\label{NB6} A linear connection $\nabla$ on $N$ defines an
effective  reductive Cartan geometry $\xi$ on $N$ of type $(G,H)$,
where $H = \GL(q,\R)$ and $G$ is the semi-direct product of the
Lie groups $\GL(q,\R)$ and $\R^q$. The Lie group $G$ is
interpreted as the Lie group $\Aff(\R^q)$ of all affine
transformations of the space $\R^q$, and $H$ is its stationary
subgroup. Thus a foliation $(M, F)$ with a transversal linear
connection is a reductive Cartan foliation.
\end{Rem}

\subsection{A linear connection associated with a reductive Cartan geometry}
Let $M$ be a smooth manifold and let $V$ be a vector space. A map
$\sigma: M \to V$ is called a $V$-valued function on $M$. Let
$\mathfrak F_V(M)$ denote the space of all $V$-valued smooth
functions on  $M$. Denote by $\sigma_{*x}$ the differential of the
map $\sigma$ at a point $x\in M$. The action of a vector field
 $X\in\mathfrak X(M)$ on  $\mathfrak
F_V(M)$ is defined by the equality
$$(X \sigma)(x)= \sigma_{*x} (X_x) \,\,\,\,\,
\forall\sigma\in\mathfrak F_V(M), \forall x\in M.$$ The map
$$X: \mathfrak F_V(M) \to \mathfrak F_V(M),\quad \sigma\mapsto
X\sigma$$ is $\R$-linear.

Let $(M, F)$ be a reductive Cartan foliation modelled on a
reductive Cartan geometry $\xi = (P(N,H), \omega)$ of type
$(G,H)$, where $\xi$ is  reductive with respect to the
decomposition of $\mathfrak g$ into the direct sum of vector
spaces $\mathfrak g = \mathfrak h \oplus V,$ here $\mathfrak g$
and $\mathfrak h$ are the Lie algebras  of the Lie groups $G$ and
$H$,  and $V$ is an $\Ad_G(H)$-invariant vector subspace of
$\mathfrak g$. The effectivity of the Cartan geometry $\xi$ is not
assumed.

Consider the smooth $q$-dimensional distribution  $$Q = \{Q_u=
\omega_u^{-1}(V)\,|\, u\in P\}$$ on $P$. The $\Ad_G(H)$-invariance
of the vector subspace $V$ of $\mathfrak g$ implies the
$H$-invariance of the distribution $Q$. Thus $Q$ is an
$H$-connection on the principal $H$-bundle $P(N,H)$. Let $p: P\to
N$ be the projection.

Recall that a $\mathfrak g$-valued function $h\in\mathfrak F(P)$
is called $H$-equivariant if it satisfies the equality $$h(ua) =
\Ad(a^{-1})h(u) \quad \forall u\in P,\,\, \forall a\in H.$$
 Denote by
$\mathfrak F_{H,V}(P)$ the set of $H$-equivariant  $V$-valued
functions on $P$. Let $\mathfrak X_{H,Q}(P)$ be the set of all
$H$-invariant smooth vector fields on $P$ tangent to $Q$. Note
that $\mathfrak F_{H,V}(P)$ and $\mathfrak X_{H,Q}(P)$ are modules
over the algebra of functions $\mathfrak F(P)$. For each vector
field $X\in\mathfrak X(M)$ there exists a unique vector field
$Y\in\mathfrak X_{H,Q}(P)$ such that $p_*(Y) = X$. We denote this
vector field by  $\widehat X$. If $h\in\mathfrak F(M)$, then let
$\widehat h= h \circ p\in\mathfrak F(P)$

Consider the map
\begin{equation}\label{eqA5}
\alpha: \mathfrak X_{H,Q}(P) \to\mathfrak F_{H,V}(P), \quad Y
\mapsto \omega(Y)\,\,\,\,\,\, \forall Y\in \mathfrak X_{H,Q}(P),
\end{equation}
which is an isomorphism of the  modules over the algebra of smooth
functions $\mathfrak F(P)$. Every vector field $Y\in\mathfrak
X_{H,Q}(P)$ defines the $\mathbb{R}$-linear map
\begin{equation}\label{eqA6}
D_Y: \mathfrak F_{H,V}(P)\to\mathfrak F_{H,V}(P) ,\quad \sigma
\mapsto Y\sigma\,\,\,\,\,\, \forall\sigma\in\mathfrak F_{H,V}(P),
\end{equation}
which satisfies the equality
\begin{equation}\label{eqA7}
D_Y (f\sigma) = (Y f)\sigma  + f D_Y (\sigma) \,\,\,\,\,\, \forall\sigma\in\mathfrak F_{H,V}(P),
\forall f\in\mathfrak{F}(P).
\end{equation}
We obtain
\begin{Prop}\label{p4} Let $\xi = (P(N,H), \omega)$ be a reductive Cartan geometry
with the projection $p: P\to N$, and $\alpha:
\mathfrak{X}_{H,Q}(P)\to\mathfrak{F}_{H,V}(P)$ be the isomorphism
 defined by $(\ref{eqA5})$. Then the equality
\begin{equation}\label{eqA8}
\nabla_X Y= p_{*}(\alpha^{-1} (D_{\widehat X} \alpha(\widehat Y)))
\,\,\,\,\,\, \forall X, Y \in \mathfrak X(N)
\end{equation}
defines a linear connection $\nabla$ on the manifold $N$.
\end{Prop}

\begin{Def} The linear connection $\nabla$ defined in Proposition~\ref{p4}
 by a reductive Cartan geometry $\xi = (P(N,H), \omega)$  is called
the linear connection associated with $\xi$.
\end{Def}

\begin{Rem}\label{NB2} Lotta \cite{Lot} proved that the  statement by Sharpe \cite[Lem. 6.4]{Shar}
about the existence of a mutation for any reductive Cartan
geometry $\xi$ of type ${\mathfrak g}/{\mathfrak h}$ with the
decomposition $\mathfrak g = \mathfrak h \oplus V$ to a Cartan
geometry $\xi'$ of type ${\mathfrak g}'/{\mathfrak h}$ with the
decomposition $\mathfrak g' = \mathfrak h \oplus \mathfrak p,$
where $\mathfrak p$ is a subalgebra in the Lie algebra ${\mathfrak
g}'$ such that $[\mathfrak p,\mathfrak p] = 0$, does not hold true
generally. Consequently the structure of the reductive Cartan
geometries is complicated than it is stated in~\cite{Shar}.
According to Proposition~\ref{p4} and Remark~\ref{NB6}, any
reductive Cartan geometry of type  $\frak g/\frak h$ with the
decomposition $\mathfrak g = \mathfrak h \oplus V$ induces a
reductive Cartan geometry  $\eta$ of type $\frak g'/\frak h'$ with
the decomposition $\mathfrak g' = \mathfrak h' \oplus\frak p$,
where $\frak p = \R^q$, $q = \dim(V) = \dim(\frak p)$, is a
subalgebra of the Lie algebra ${\mathfrak g}'$, which is the Lie
algebra of $\GL(q,\R)$, hence $[\mathfrak p,\mathfrak p] = 0$. In
general $\mathfrak h\neq\mathfrak h'$ and $\eta$ is not a mutation
of $\xi$.
\end{Rem}

\subsection{A transversal linear connection associated to a
reductive Cartan foliation}

\begin{Th}\label{th6} Each reductive Cartan foliation is a foliation with a transversal linear connection.
\end{Th}

\begin{Cor}\label{c2} The holonomy group of each leaf of a reductive Cartan foliation is linearizable.
\end{Cor}

\begin{Rem} \label{NB1} The topology of a foliated manifold $M$ is invariant under the change of the
transversal geometric structure of the  foliation $(M, F)$ by any
other geometric structure on this foliation. Hence, due to
Theorem~\ref{th6}, any topological problem for a reductive
Cartan foliation may be reduced to a similar problem for a
foliation with a transversal linear connection.
\end{Rem}

{\bf Proof of Theorem~\ref{th6}.}

\begin{Lem} \label{l1} Let $\xi = (P(N,H), \omega)$ and $\xi' = (P'(N',H), \omega')$ be
two reductive Cartan geometries of the same type $(G,H)$ with
respect to the decomposition
 $\frak g  = \frak h \oplus V$, having the associated  linear connections $\nabla$ and $\nabla'$,
respectively. If $\Gamma: P\to  P'$ is an isomorphism of the
Cartan geometries $\xi$ and $\xi'$, then its projection $\gamma:
N\to N'$ is an isomorphism of the associated  linear connections
 $\nabla$ and $\nabla'$.
\end{Lem}

{\it Proof.} All objects related to the geometry $\xi'$ will be
denoted with the prime. Since  $\Gamma: P\to  P'$ is an
isomorphism of the Cartan geometries $\xi$ and $\xi'$, it holds
$\Gamma^*\omega' = \omega.$ Hence, $\Gamma_* Q = Q'$. Consequently
the isomorphisms of the vector spaces
$$\Gamma^* :\mathfrak F_{H,V}(P') \to \mathfrak F_{H,V}(P),\quad \sigma \mapsto \sigma\circ\Gamma,
\quad \Gamma_* :\mathfrak X_{H,Q}(P) \to \mathfrak
X_{H,Q'}(P'),\quad Y\mapsto \Gamma_* Y$$ are induced. These
isomorphisms satisfy the equality $$(\Gamma^*)^{-1}\circ\alpha =
\alpha'\circ\Gamma_*,$$ where $\alpha$ and $\alpha'$ were defined
above. Therefore we get
$$\widehat{\gamma_*\nabla_X Y} = \Gamma_* \alpha^{-1}(D_{\widehat
X} \alpha(\widehat Y)) =\alpha^{' {-1}}\circ(\Gamma^*)^{-1}(D_{\widehat
X} \alpha(\widehat Y)) = \alpha^{' {-1}} (D _{\Gamma_*\widehat X}
\alpha'(\Gamma_*\widehat Y))
 =\widehat{\nabla '_{\gamma_*X}\gamma_*Y},$$ that means that,
$$\widehat{\gamma_*\nabla_X Y} = \widehat{\nabla
'_{\gamma_*X}\gamma_*Y}.$$ Since the map  $\mathfrak X(N) \to\frak
X_{H,Q}(P): X \mapsto\widehat {X}$ is bijective, it holds
$$\gamma_*\nabla_X Y = \nabla '_{\gamma_*X}\gamma_*Y$$ for all  $X,
Y \in\mathfrak X(N)$. Thus $\gamma_*: N\to N'$ is an isomorphism
of the manifolds with linear connections $(N,\nabla)$ and
$(N',\nabla')$. \qed

Let $(M, F)$ be a reductive Cartan foliation modelled on a
reductive Cartan geometry $\xi = (P(N,H), \omega)$. Suppose that
$(M, F)$ is defined by an $(N,\xi)$-cocycle $\{U_i,f_i,
\Gamma_{ij}\}_{i,j \in J}$. According to Proposition~\ref{p4}, the
associated linear connection $\nabla$ is defined on the manifold
$N$. Since every transformation $\Gamma_{ij}$ is a local
isomorphism of the corresponding reduced Cartan geometries, by
Lemma~\ref{l1}, its projection $\gamma_{ij}$ is an isomorphism of
the induced linear connections on the open subsets of $N$. This
means that $(M, F)$ is a foliation with  a transversal linear
connection given by the $(N,\nabla)$-cocycle $\{U_i, f_i,
\gamma_{ij}\}_{i,j \in J}$. This proves Theorem~\ref{th6}. \qed

\section{Attractors and minimal sets}\label{secattrminsets}

\subsection{Statement of the results}
Recall that {\it a minimal set} of a foliation on a manifold $M$
is a nonempty closed subset in $M$ that consists of a union of
leaves and has no proper subset satisfying this condition. Minimal
sets for transformation groups are defined in a similar way.

\begin{Rem}\label{NB3} A nonempty closed saturated subset $\mathcal M$ in $M$ is an
attractor and a minimal set of the foliation $(M, F)$ if  and only
if there exists an open neighbourhood  $\mathcal U = \mathcal
U(\mathcal M)$ of the subset $\mathcal M$ such that
$\overline{L}\supset\mathcal M$ for any leaf $L\subset\mathcal U.$
\end{Rem}

\begin{Def}\label{d3} Let $\Phi$ be a  group of homeomorphisms of a topological space
$N$. We call a point $x\in N$ a local limit point of the  group
$\Phi$ if there exists a neighbourhood $U$ of the point $x$ such
that the closure of the orbit $\Phi{.}y$ of any point $y\in U$,
$y\notin\Phi{.}x$, contains $x$. If, moreover, $U = N$, then  $x$
 is called the limit point of the group $\Phi$.
\end{Def}

\begin{Rem} \label{NB4} The origin $0\in\mathbb R^q$ is a fixed point
for any subgroup $\Phi$ of the linear group $\GL(q,\R)$.
Consequently the origin is the only possible local limit point of
$\Psi$.
\end{Rem}

Consider a reductive Cartan foliation $(M, F)$ of codimension $q$
and a transversal $q$-dimen\-sio\-nal distribution $\mathfrak M$.
Let ${\mathcal R}(M,H)$ be the principal $H$-bundle with the
$H$-invariant connection $Q$, defined by above the transversal
reductive Cartan geometry $\xi$  and let ${\mathcal
R'}(M,\GL(q,\R))$ be the principal $\GL(q,\R)$-bundle with the
$\GL(q,\R)$-invariant connection $Q'$ given by the transversal
geometry $(N,\nabla)$. Note that an integral curve $\sigma$ of the
distribution $\mathfrak M$ is a geodesic of the connection $Q$ if
and only if $\sigma$ is a geodesic of the connection $Q'$. Such
curves are called $\mathfrak M$-geodesics.  We observe that the
completeness of a reductive Cartan foliation is equivalent to the
existence of a transversal distribution $\mathfrak M$ such that
every maximal $\mathfrak M$-geodesic is defined on the whole real
line.

We give now sufficient conditions for the existence of an
attractor (and a global attractor) that is also a minimal set.

\begin{Th} \label{th7} Let $(M, F)$ be a reductive Cartan foliation of codimension
$q$. Suppose that there exists a leaf $L$ such that its linear
holonomy group $D\Gamma(L,x)$ at some point $x\in M$ has a limit
point. Then:
\begin{itemize} \item[(1)] The closure of the leaf ${\mathcal M} = \overline{L}$
is an attractor and a minimal set of the foliation $(M, F)$.
\item[(2)] If, moreover, $(M, F)$ is a complete Cartan foliation
with respect to a transversal $q$-dimen\-sio\-nal distribution
$\mathfrak M$ and the leaf $L$ can be connected with every leaf
$L_\alpha$ of $(M, F)$ by a smooth $\mathfrak M$-geodesic, then
${\mathcal M} = \overline{L}$ is a global attractor and a minimal
set of this foliation. \end{itemize}
\end{Th}

\begin{Cor}\label{c3}  Let $(M, F)$ be a complete reductive Cartan foliation of codimension $q$.
Suppose that the curvature and the torsion of the associated
transversal linear connection are zero. If there exists a leaf $L$ with
a linear holonomy group admitting a limit point,
then:
\begin{itemize} \item[(1)] the closure of the leaf ${\mathcal M} = \overline{L}$
is a global attractor and a minimal set of the foliation $(M, F)$;
\item[(2)] there exists a regular covering map $\kappa: \widehat{M}\to M$ such that
the foliation $(M, F)$ is covered  by the trivial bundle $r:
\widehat{M} = L_0\times \R^q\to \R^q$ over the  space $\R^q$,
where $L_0$ is a manifold diffeomorphic to every leaf without
holonomy;
\item[(3)] the global holonomy group $\Psi$ of $(M, F)$ is a subgroup of the affine Lie group
$\Aff(\R^q)$, it has a global attractor $\mathcal K$, and
${\mathcal M} = \kappa(r^{-1}(\mathcal K)).$
\end{itemize}\end{Cor}

\begin{Cor}\label{c4}  Let $(M, F)$ be a  reductive Cartan foliation of codimension $q$.
Suppose that there exists a leaf $L$ such that its linear holonomy
group contains an element defined by a matrix of the form $B\cdot
D$, where $B\in O(q)$ and $D = {\rm diag} (d_1,...,d_q)$ with
$|d_i| < 1$ for $1\leq i\leq q$. Then $(M, F)$ has an attractor
${\mathcal M} = \overline{L}$ which is a minimal set.
\end{Cor}

Recall that a smooth foliation  $(M, F)$ is called {\it proper} if
each its leaf is an embedded submanifold in $M$. A leaf $L$ is
called closed if it is a closed subset in $M.$ As it is known,
(see e.g. \cite{Tam}), any minimal set of a foliation is either a
closed leaf, or the closure of a non-proper leaf. This and
Theorem~\ref{th7} imply the following statements:

\begin{Cor} \label{c5} Let $(M, F)$ be a reductive Cartan foliation.
If there exists a proper leaf $L$ with a linear holonomy group
admitting a limit point, then the leaf $L$ is  closed and it is an
attractor of the foliation $(M, F)$.
\end{Cor}

\begin{Cor}\label{c6} Let $(M, F)$ be a complete reductive Cartan foliation.
Suppose that the curvature and the torsion of the associated
transversal linear connection are zero. If there exists a proper
leaf $L$ with a  linear holonomy group admitting a limit point,
then $L$ is a unique closed leaf, and $L$ is a global attractor of
the foliation.
\end {Cor}

\subsection{The existence of an attractor which is a minimal set}
Denote by  $\Gamma (L,x)$ the germ holonomy group of an arbitrary
leaf $L=L(x)$ at a point $x$; $\Gamma (L,x)$ consists of germs of certain
holomorphic diffeomorphisms  $\psi$ of a transversal
$q$-dimensional disk $D_x^q$ at the point \cite{Tam}. Let $D\Gamma
(L,x)$ be the linear holonomy group consisting of the
differentials  $\psi_{*x}: \mathfrak{M}_x\to\mathfrak{M}_x$, where
$\mathfrak{M}_x = T_xD_x^q.$

Suppose that the linear holonomy group $D\Gamma (L,x)$ of a leaf
$L= L(x)$ has a limit point. There exists a submersion  $f: U\to
V$ from an $(N,\nabla)$-cocycle $\{U_i, f_i, \gamma_{ij}\}_{i,j
\in J}$ defining the foliation $(M, F)$ such that $x\in U.$ Let $v
= f(x)\in V.$ Denote by  $\mathcal H$ the holonomy pseudogroup
generated by the local automorphisms $\gamma_{ij}$, $i, j\in J$ of
the transversal manifold with the linear connection $(N,\nabla)$.
Let $${\mathcal H}_v= \{\phi\in{\mathcal H}\,|\,\phi(v) = v\}.$$
We consider the linear holonomy group $D\Gamma (L,x)$ of the leaf
$L$ as the group of linear transfor\-ma\-tions $${\mathcal
H}_{*v}= \{\phi_{*v} \,| \, \phi\in{\mathcal H}_v\}$$ of the
tangent space $T_vN$ at the point $v\in V\subset N.$

It is well known that the group $\Gamma (L,x)$ is isomorphic to
the group of germs of local automorphisms $\phi\in {\mathcal H}_v$
at the point $v$. Since a linear connection defines a
$G$-structure of the first order, there exists an isomorphism
$$d_v: \Gamma (L,x) \to D\Gamma (L,x)\cong {\mathcal H}_{*v}:
\{\phi\}_v \mapsto\phi_{*v}\quad \forall\phi\in {\mathcal H}_v,$$
assigning to a germ $\{\phi\}_v $ at the point $v$ the
differential $\phi_{*v}$.

Let $W_0$ be a normal neighbourhood of the origin in the tangent
space  $T_vN$. The exponen\-ti\-al map
$$\Exp_v: W_0 \to W: X \mapsto \gamma_X(1)$$ defined by the  linear
connection $\nabla$ is a diffeomorphism onto an open neighbourhood
$W$  of $v$ in $N$. By the property of the exponential map, each
transformation $\phi\in{\mathcal H}_v$ satisfies {\bf the} following
equality
\begin{equation}\label{eq9}
\phi\circ \Exp_v =  \Exp_v  \circ\phi_{*v}
\end{equation}
in a neighbourhood of the origin in  $T_vN$, where the both sides
of Equality \eqref{eq9} are defined. Note that the holonomy
pseudogroup $\mathcal H$ of a Cartan foliation is
quasi-analytical, i.e., if a transformation $\phi\in{\mathcal
H}_v$ equals  the identity on an open subset in  $N$, then it
coincides with the identity transformation everywhere in the
domain of its definition. Since the differential $\phi_{*v}$ of
each $\phi\in{\mathcal H}_v$ is defined on the whole tangent space
$T_vN$, Equality \eqref{eq9} allows us to extend each local
automorphism $\phi\in{\mathcal H}_v$ to the whole neighbourhood
$W$ of the point $v$. Since $\mathcal H$ is quasi-analytical, this
extension is defined uniquely. Thus we assume that
Equality~\eqref{eq9} holds on~$W_0$.

By the assumption, the group ${\mathcal H}_{*v}$ has a limit
point, hence this point is the  origin in $T_{v}N$. Therefore for
any vector $Y\in T_vN$ there exists a sequence $\phi_k \in
{\mathcal H}_{v}$ such that $(\phi_{k})_{*v} (Y) \to 0$ as $k \to
+ \infty.$ Without loss generality we assume that
$(\phi_{k})_{*v}(Y) \in W_0$ for any $Y\in W_0.$  We introduce the
notation $U_0= f^{-1}(W)$ and
$${\mathcal U} = \bigcup_{ L_\alpha\in F, \, L_\alpha\cap U_0 \neq\emptyset}  L_\alpha.$$

Consider any leaf $L_\alpha\subset{\mathcal U}$. Let $x_\alpha\in
L_\alpha\cap U_0$. Then $y= f(x_\alpha)\in W.$ Since the map
$\Exp_v|_{W_0}: W_0\to W$ is a diffeomorphism, for each $y\in W$
there exists a vector  $Y = \Exp_v ^{-1}(y)\in W_0$. From
(\ref{eq9}) it follows that  $\phi_k(y) \to v$ as $k \to +
\infty.$ We emphasize that the set $f^{-1}(\{\phi_k(y)\,|\,
k\in\mathbb N\})$ is contained in $L_\alpha\cap U.$ Consequently,
$L\subset\overline{L_\alpha}$. Thus it holds
$${\mathcal M}=
\overline{L}\subset\overline{L_\alpha}\quad \forall
L_\alpha\subset \mathcal U,$$ i.e. ${\mathcal M}$ is an attractor
with the basin ${\mathcal U}$.

Let us show that ${\mathcal M}$ is a minimal set of the foliation
$(M, F)$. Let $L_\alpha$ be any leaf of the  foliation contained
in ${\mathcal M}$, i.e., $L_\alpha\subset\mathcal M$. Then the
closure  $\overline{L_\alpha}$ satisfies
$\overline{L_\alpha}\subset\overline{L} = \mathcal M.$ Since
$L_\alpha\subset\mathcal U$, from the above it follows that
$\overline{L_\alpha}\supset\overline{L}$. Consequently,
$\overline{L_\alpha}=\overline{L} = \mathcal M.$ This means that
${\mathcal M}$ is a minimal set of the foliation $(M, F)$. Thus
the statement $(1)$ of Theorem \ref{th7} is proved.

\subsection{The existence of a global attractor which is a minimal set}

We use the notations from Section~\ref{ss3.2A}. Let $(M, F)$ be a
$\mathfrak M$-complete reductive Cartan foliation. Let
$\widetilde{\mathfrak M}=\pi^*\mathfrak M,$ where $\pi: {\mathcal
R}\to M$ is the associated $H$-bundle. Let $g_{\mathcal R}$~be a
Riemannian metric on $\mathcal R$. Consider an Euclidean metric
$d_0$ on the vector space $\mathfrak g$ such that $\mathfrak h$
and $V$ are orthogonal subspaces. Let $Z=Z_{\mathcal F}\bigoplus
Z_{\widetilde{\mathfrak M}}$ be the decomposition of a vector
field $Z\in \mathfrak X({\mathcal R})$ with respect to the
decomposition $T_u{\mathcal R}= T_u{\mathcal F}\oplus
\widetilde{\mathfrak M}_u, u\in \mathcal R.$ The equality
$$d(X,Y)= g_{\cal R}(X_{\mathcal F},Y_{\mathcal F})+d_0({\omega}(X),
{\omega}(Y)), \forall X,Y\in \mathfrak X({\cal R}),$$ defines a
Riemannian metric $d$ on $\mathcal R,$ and $d$ is transversally
projectible with respect to the lifted foliation $({\mathcal
R},\mathcal F)$.

Let $E_{i}, i=1,\dots,\dim\mathfrak g$ be a basis of the Lie
algebra $\mathfrak g$. Denote by $ X_{i}$ the vector field from
$\mathfrak{X}_{\widetilde{\frak M}}(\mathcal R)$ such that
$\omega(X_i) = E_i.$  Let $\widehat{\nabla}$ be the Levi-Civita
connection of the Riemannian manifold $({\mathcal R},d).$ It can
be checked  directly that the equality
\begin{equation}\label{eq*} \widetilde{\nabla}_{Y}{Z}=Y({Z}^{i})X_{i}+\widehat{\nabla}_{Y}{Z}_{\mathcal
F},\end{equation} where $Z=Z_{\mathcal F}\oplus Z_{\widetilde
{\mathfrak M}}$, $Z_{\widetilde{\mathfrak
M}}={Z}^{i}X_{i}\in\mathfrak X_{\widetilde{\mathfrak M}}({\mathcal
R})$, $Z_{\mathcal F}\in \mathfrak X_{\mathcal F}({\mathcal R})$,
and $Y\in\mathfrak X(\mathcal R)$, defines a linear connection
$\widetilde{\nabla}$ in $\mathcal R$. The connection
$\widetilde{\nabla}$ is in general not torsion free,  and it holds
$\widetilde{\nabla} d = 0$. Every vector field
$X\in{\mathfrak{X}}_{\widetilde{\mathfrak M}}(\mathcal{R})$ such
that $\omega(X)=const\in\mathfrak g$, is parallel with respect to
$\widetilde{\nabla}$, hence its integral curves are geodesic lines
of $\widetilde{\nabla}$.

Let $$Q = \{Q_u=\omega_u^{-1}(V)\,|\, u\in\mathcal R\}$$ be the
$H$-connection in $\mathcal R$ and let $$\mathfrak N = \{\mathfrak
N_u= Q_u\cap\widetilde{\mathfrak M}_u\,|\, u\in\mathcal R\}.$$
Geodesics of $\widetilde{\nabla}$ that are integral curves of the
distribution $\mathfrak N$ are called $\mathfrak N$-geodesics.
Note that for any $\mathfrak N$-geodesic $\gamma$, the curve
$\pi\circ\gamma$ is an $\mathfrak M$-geodesic.

As the Cartan foliation $(M, F)$ is $\mathfrak M$-complete, the
exponential map $\Exp_u$, $u\in\mathcal R$, of
$\widetilde{\nabla}$ is defined, in particular, on $\mathfrak
N_u$. Therefore the map $\Exp_x: \mathfrak M_x\to M$, $x=\pi(u)\in
M$, satisfying the equality $\Exp_x\circ\pi_{*u} = \pi\circ
\Exp_{u}$, is defined.

Since $(M, F)$ admits a leaf $L = L(x_0)$ such that its linear
holonomy group $D\Gamma(L,x_0)$ has a limit point, according to
the proved statement $(1)$ of Theorem~\ref{th7}, the closure
${\mathcal M} =\overline{L}$ is an attractor. Let ${\mathcal B} =
{\mathcal B}({\mathcal M})$ be its basin. Then there exists an
open star neighbourhood $V_0$ of zero  in ${\mathfrak M_{x_0}}$
such that $\Exp_{x_0}(V_0)\subset{\mathcal B}$.

By the assumptions, for any leaf $L_\alpha$, there exists an
$\mathfrak M$-geodesic $\sigma: [0,1]\to M$ such that $x_0 =
\sigma(0)\in L$ and $x = \sigma(1)\in L_\alpha$. Pick
$u_0\in\pi^{-1}(x_0)$, then there exists an $\mathfrak N$-geodesic
$\gamma: [0,1]\to\mathcal{R}$ which is a $\mathfrak N$-lift of
$\sigma$ starting at the point $u_0$, i.e. $\gamma(0) = u_0$ and
$\sigma = \pi\circ\gamma$. The $\mathfrak M$-completeness of $(M,
F)$ implies the existence of a vector $Y\in\mathfrak N_{u_0}$ such
that $\gamma(s) = \Exp_{u_0}(sY)$ for any $s\in [0,1]$. Therefore
the vector $X = \pi_{u_0}(Y)$ satisfies the relation $\sigma(s) =
\Exp_{x_0}(sY)$ for any $s\in [0,1]$.

Since the linear holonomy group $D\Gamma(L,x_0)$ has a limit
point, there exists an element $\phi_{*x_0}\in D\Gamma(L,x_0)$ for
which  $Y= \varphi_{*x_0}(X)\in V_0$. There is a loop $h: [0,1]\to
L$ at the point $x_0$ such that $\varphi$ is a local holonomic
diffeomorphism of a transversal $q$-dimension disk $D_{x_0} =
\Exp_{x_0}(V_0)$ along $h$, and $\varphi(x_0) = x_0$. Let
${\mathcal L} = {\mathcal L} (u_0)$ be a leaf of the lifted
foliation $(\mathcal R,\mathcal F)$. Since $\pi_{\mathcal L} :
{\mathcal L}\to L$ is a covering map, there exists a curve
$\widetilde{h}: [0,1]\to{\mathcal L}$ starting at
$\widetilde{h}(0) = u_0$ and covering $h$, i.e.,
$\pi\circ\widetilde{h} = h.$

Recall that $\widetilde{\mathfrak M}$ is an Ehresmann connection
for the foliation $(\mathcal R,\mathcal F)$. Let $\gamma\stackrel
{\widetilde{h}}{\mapsto}{\gamma}^*$ and $\widetilde{h}\stackrel
{\gamma}{\mapsto}\widetilde{h}^*$ be the translations with respect
to the Ehresmann connection $\widetilde{\mathfrak M}$. Consider
the translations $\sigma\stackrel {h}{\mapsto}{\sigma}^*$ and
$h\stackrel {\sigma}{\mapsto}h^*$ with respect to the Ehresmann
connection ${\mathfrak M}$. Then $\widetilde{h}^*$ is a curve in a
leaf $\widetilde{\mathcal L}$ of $(\mathcal R,\mathcal F)$ and
$h^*$ is a curve in a leaf $L_\alpha,$ with
$\pi(\widetilde{\mathcal L}) = L_\alpha$ and
$\pi\circ\widetilde{h}^* = h^*$. Note that $\sigma^*(s) =
\Exp_{x_0}(sY)\in\mathcal{B}$  for all $s\in [0,1]$. Therefore, $y
= h^*(1) = \sigma^*(1)\in L_\alpha\cap\mathcal{B}$. Since
$\mathcal{B}$ is a saturated set, we have
$L_\alpha\subset\mathcal{B}$. Thus, taking into attention the fact
that $L_\alpha$ is an arbitrary leaf of $(M, F)$, we get $M =
\mathcal{B}$, i.e. $\mathcal M$ is a global attractor. This
completes the proof of the statement $(2)$ of Theorem \ref{th7}.
Theorem \ref{th7} is proved. \qed

\paragraph{Proof of Corollary~\ref{c3}.}
Since the curvature and the torsion of the linear connection
$\nabla$ are equal to zero, $(N,\nabla)$ is a locally affine
manifold. Therefore $(M, F)$ is a transversally affine foliation.
In other words, $(M, F)$ is an $(\Aff(\R^q),\R^q)$-foliation.
Applying Theorem \ref{th4} to the complete
$(\Aff(\R^q),\R^q)$-foliation $(M, F)$, we see that there exists a
regular covering map $\kappa:\widehat{M}\to M$ such that the
leaves of the induced foliation $(\widehat{M},\widehat F)$ are
fibres of a locally trivial bundle $r:\widehat{M}\to B$ over a
simply connected affine manifold $B$.

Let $\mathfrak M$ be a transversal distribution on $M$ with
respect to which the reductive Cartan foliation $(M, F)$ is
$\mathfrak M$-complete. Then $\mathfrak M$ is an Ehresmann
connection for $(M, F)$. This implies that $\widehat{\mathfrak M}=
\kappa^*\mathfrak M$ is an Ehresmann connection for the foliation
$(\widehat{M},\widehat F)$. In this case $\widehat{\mathfrak M} $
is an Ehresmann connection for the submersion $r:\widehat{M}\to M$
\cite[Prop. 2]{BH}. Hence any geodesic on $B$ admits an
$\widehat{\mathfrak M}$-lifts in $\widehat{M}$. Since every such
lift $\widehat{\sigma}$ of a maximal geodesic $\sigma$ from $B$ is
a maximal $\widehat{\mathfrak M}$-geodesic in $\widehat{M}$, the
canonical parameter on $\sigma$ is defined on the real line. This
means that the  affine manifold $B$ is complete. Thus $B$ is a
simply connected complete torsion free affine manifold with zero
curvature tensor. Therefore $B$ is the  affine space $\R^q$. Every
locally trivial bundle over a contractible base is trivial, hence
we get the trivial bundle $r: \widehat{M} = L_0\times \R^q\to
\R^q$, and the manifold $L_0$ is diffeomorphic to any leaf of $(M,
F)$ with the trivial holonomy group.

Assume that there exists a leaf $L$ of $(M, F)$ such that its
linear holonomy group has a limit point. According the statement
$(1)$ of Theorem~\ref{th7}, the closure ${\mathcal M} =
\overline{L}$ is an attractor and a minimal set of  $(M, F)$.

Since any two points in $\R^q$ may be connected by a geodesic,
using the Ehresmann connection $\widehat{\mathfrak M}$ we able to
connect any two leaves of $(\widehat{M},\widehat{F})$ by an
$\widehat{\mathfrak M}$-geodesic. Therefore for every leaf
$L_\alpha$ of $(M, F)$ there exists a ${\mathfrak M}$-geodesic
connecting $L$ with $L_\alpha$. From this and the statement $(2)$
of Theorem~\ref{th7} it follows  that ${\mathcal M}$ is a global
attractor of the foliation $(M, F)$. \qed

\section{Examples}\label{s7}

First of all we stress that the foliations admitting transversally
projectible Riemannian metrics do not admit attractors \cite{Min}.
As it is known \cite[Th. 4]{Min}, Cartan foliations of type
$\mathfrak g/\mathfrak h$, where the Lie algebra $\mathfrak h$ is
compactly embedded into the Lie algebra $\mathfrak g$, are
Riemannian foliations, hence they also do not admit attractors.
Example~1 below provides a complete transversely affine foliation
that does not admit an attractor and that is not a Riemannian
foliation. Example~2 shows that in the framework of
Theorem~\ref{th7}, the situation ${\mathcal M} = \overline{L} = M$
is possible. According to Example~3, there exist non-complete
transversal affine foliations with global attractors. In Example~4
we construct a transversally affine foliation with regular global
attractor that illustrates Corollary~\ref{c3}. Examples of global
attractors of transversally similar foliations are constructed
in~\cite{Min}.

 We denote by $f_A = \langle A,a\rangle $ the element of the affine group
$\GL(q,\R)\ltimes \R^q\cong \Aff(\R^q)$, where $A\in \GL(q,\R)$
and $a\in \R^q$. It holds  $\langle A,a\rangle x = Ax + a$ for any
$x\in \R^q$ and $\langle A,a\rangle \circ\langle B,b\rangle  =
\langle AB, Ab+a\rangle $ for the composition of every two
elements $\langle A,a\rangle , \langle B,b\rangle \in\Aff(\R^q)$.

\subsection*{Example 1}
Let $f_A: \R^2\to \R^2$ be the affine transformation of the plain
$\R^2$ given by the matrix $A$=$\left(
                                         \begin{array}{ccc}
                                          1/2 & 0 \\
                                          0 & 2\\
                                         \end{array}\right)$
with respect to the canonical basis $e_1 = \left(
                                         \begin{array}{cc}
                                          1 \\ 0 \\
                                         \end{array}\right)$,
$e_2 = \left(
                                         \begin{array}{cc}
                                          0 \\ 1 \\
                                         \end{array}\right)$.
Let $B = \S^1$ be the unite circle. We define the group
homomorphism
 $\rho: \pi_1(\S^1,b)\to \Aff(\R^2)$
by setting its value on the generator
$1\in\mathbb{Z}\cong\pi_1(\S^1,b)$ to $\rho(1) = f_A.$ Then the
suspended  foliation $(M, F)= \Sus(\R^2,B,\rho)$ is a transversely
affine foliation. This foliation is $\mathfrak M$-complete, where
$\mathfrak M$ is the tangent distribution to the transversal
locally trivial bundle $p:M\to \S^1.$

The foliation $(M, F)$ is covered by the trivial bundle
$r:\R^1\times \R^2\to \R^2$. It is easy to see that its  global
holonomy group equals to $\Psi = <f_A>\cong\mathbb{Z}$. For each
 $v\in M$ there exists a point $z\in r(p^{-1}(v))\in \R^2$.
Note that the leaf  $L=L(v)$ is compact and it is diffeomorphic to
the circle if an only if $z\in\{(x,y)\in \R^2 \,|\, xy\neq 0\}
\cup \{(0,0)\}$. Consequently, the  global holonomy group $\Psi$
has no attractors. A foliation $(M, F)$ defined as above by a
matrix $A$
 is Riemannian if and only if $A$ belongs to the orthogonal group
$O(2)$. In our case $A\notin O(2)$. Thus the complete transversely
affine  foliation $(M, F)$ does not admit an attractor and it is
not a Riemannian foliation.

\subsection*{Example 2.} Let, as above,  $e_1$, $e_2$ be the canonical basis of the plain $\R^2$.
Let us consider real numbers $\lambda_0, \lambda_1, \lambda_2$
such that $0< \lambda_i < 1, \,i=0, 1, 2.$ Let $f_0(x)=
\lambda_0(x)$, $f_1(x)= \lambda_1(x - e_1) +e_1$, $f_2(x)=
\lambda_2(x - e_2) +e_2 $ $\,\, \forall x\in \R^2$.

Denote by $B = \mathbb S^2_3$ the sphere with three handles. For
its fundamental group we have
$$\pi_1(B) = < a_j,b_j\,|\, j= 0,1,2;\,\, a_0b_0a_0^{-1}b_0^{-1}a_1b_1a_1^{-1}b_1^{-1}a_2b_2a_2^{-1}b_2^{-1}=1>.$$
Define the group homomorphism $\rho: \pi_1(B)\to \Diff (\R^2)$
assuming that
 $$\rho(b_j) = \id_{\R^2}, \quad \rho(a_j) = f_{j}, \,\, j= 0, 1, 2.$$
We get the suspended  foliation $(M, F) = \Sus (\R^2,\mathbb
S^2_3,\rho)$ on the noncompact $4$-dimensional manifold $M$, which
is the total  space of a  locally trivial bundle over $B$ with the
standard fibre~$\R^2$.

By \cite[Prop. 16]{Min}, the orbit $\Psi.x$ of any point
$x\in\R^2$ is dense in $\R^2$. Thus  $(M, F)$ is a transversely
affine foliation, and $M$ is its  minimal set.

\subsection*{Example 3.} Denote by $0_q$ the origin in $\R^q$. Suppose that $q\geq 3$.
Consider the submanifold $\widehat{M} = \R^q\setminus\{0_q\}$ of
$\R^q$. Let $(\widehat{M},\widehat{F})$ be the simple foliation
defined by the submersion  $$r: \widehat{M}\to \R^1,\quad
(x_1,...,x_q)\mapsto x_1.$$ Consider the homothety  $\gamma =
\langle \lambda E, 0\rangle $ with the coefficient $\lambda>1$,
where $E$ is the identity matrix of order $q$. We obtain the
affine Hopf manifold $M= \widehat{M}/\Gamma,$ where $\Gamma =
<\gamma>$ is a group of similarity transformations of
$\widehat{M}$; note that $M\cong \S^{q-1}\times \S^1$. Since
$r\circ\gamma = \gamma\circ r$, on the manifold  $M$ the foliation
$F$ is induced such that its leaves are images of the leaves of
the foliation $(\widehat{M}, \widehat{F})$ under the universal
covering $\kappa: \widehat{M}\to M = \widehat{M}/{\Gamma}$.
Consequently the foliation $(M, F)$ is covered by the simple
foliation $(\widehat{M}, \widehat{F})$, and both these foliations
are transversally affine. Note that the leaf $L_0 =
\kappa(r^{-1}(0_1))$ is compact and diffeomorphic to
 $\S^{q-2}\times \S^1$, and it is a  global
attractor of the  foliation $(M, F)$. All the other leaves of this
foliation are diffeomorphic to $\R^{q-1}$.

Suppose that the foliation  $(M, F)$ admits an Ehresmann
connection. Then by Proposition~\ref{p2}, it is covered by a
locally trivial bundle  $p: \widehat{M}\to \R^1$, whose fibres are
all diffeomorphic to each other. Since not all leaves of the
submersion $r: \widehat{M}\to \R^1$ are diffeomorphic to each
other, we get a contradiction. Thus   $(M, F)$  does not admit an
Ehresmann connection, consequently it is not a  complete
transversely affine foliation.

The linear holonomy group of the leaf $L_0$ has a limit point,
hence the foliation $(M, F)$ with a global attractor satisfies the conditions of
the statement $(1)$ of Theorem~\ref{th7}, but it does not satisfy the conditions of
the statement $(2)$ of Theorem~\ref{th7}.

\subsection*{Example 4.} Let $\R^2$ be the plane with the coordinates $x,\, y.$
Consider two affine transformations of $\R^2$: $c= \langle
A,0\rangle $, where $A$=$\left(\begin{array}{ccc}
                                         \mu_1 & 0 \\
                                          0 & \mu_2\\
                                         \end{array}\right)$
and $\psi_2 = \langle C,c\rangle $, where
$C$=$\left(\begin{array}{ccc}
                                         \mu_3 & 0 \\
                                          0 & \nu\\
                                         \end{array}\right)$, $\mu_i$, $i = 1, 2, 3$, $\nu$  are real numbers
                                                                                such that $|\mu_i|<1$ and   $c = \left(
                                         \begin{array}{cc}
                                          1 \\ 0 \\
                                         \end{array}\right)$. The point $x_0 = \left(
                                         \begin{array}{cc}
                                          \frac{1}{1-\mu_3} \\ 0 \\
                                         \end{array}\right)$ is the only fixed point of $\psi_2$.
Denote by $\Psi$ the subgroup of the affine group $\Aff(\R^2)$
generated by $\psi_1$ and $\psi_2$. We show that $\Psi$ has a
global attractor $\mathcal K$ coinciding with the coordinate axis
$x$. Let
$$h_n= \psi_1^{n}\circ\psi_2^{n}\circ\psi_1^{- n}\circ\psi_2^{-
n}.$$ It can be shown that $$h_n= \langle E,d_{n}\rangle,$$ where
$d_n = \left(
                                         \begin{array}{cc}
                                         \frac{(\mu_1^{n}-1)(\mu_3^{n}-1)}{\mu_3-1}\\ 0 \\
                                         \end{array}\right)$ for every $n\in\mathbb N$.
It can be checked directly that $$h_{n+1}\circ h_{ n}^{-1} =
\langle E,\delta_{n}\rangle ,$$ where $\delta_n = \left(
                                         \begin{array}{cc}
                                           \delta_n^{(1)}\\ 0 \\
                                         \end{array}\right)$, and $\delta_n^{(1)} =
\frac{\mu_1^{n+1}(\mu_3^{n+1} - 1) - \mu_1^{n}(\mu_3^{n} -
1)}{\mu_3 - 1} - \mu_3^n$. Hence, $\delta_n\to 0$ as $n\to
+\infty.$ Taking into account the fact that $h_{n+1}\circ
h_n^{-1}\in\Psi$ is a parallel translation along the axis $x$, we
see that the orbit $\Psi.0$ of the origin of $\mathbb \R^2$ is
dense in the axis $x$, and the closure of this orbit ${\mathcal
K}$ is a global attractor of the group $\Psi$.

Let $B = (\S^1\times \S^2) \sharp (\S^1\times \S^2)$ be the
connected sum of two copies of the product $\S^1\times \S^2$ of
the unit circle $\S^1$ and the unite two-sphere $\S^2$. The
fundamental group $\pi_1(B) = <g_1, g_2>$ is a free group of rank
two.

Consider the group homomorphism $\rho: \pi_1(B)\to\Psi$ defined by
the conditions $\rho(g_i)  = \psi_i$, $i = 1,2$. It defines the
suspended foliation $(M, F) = \Sus(\R^2, B,\rho)$ with the global
holonomy group $\Psi$. Let $\kappa: \widehat{M}\to M$ be the
regular covering map satisfying the conditions of
Theorem~\ref{th4}. According to Corollary~\ref{c3},
$\widehat{M}\cong L_0\times\R^2$, where $L_0$ is diffeomorphic to
every leaf without holonomy, and leaves of the induced foliation
$(\widehat{M},\widehat{F})$, $\widehat{F} = \kappa^*F$, are fibres
of the trivial bundle $r: L_0\times \R^2\to \R^2$.

Consequently ${\mathcal M}= \kappa(r^{-1}(\mathcal K))$   is a
global attractor of the complete transversally affine foliation
$(M, F)$, and $\mathcal M$ is an embedded submanifold of $M$.
Since $M$ is the space of a locally trivial bundle $p: M\to B$
with the contractible standard fibre $\mathbb R^2$, the
 exact homotopy sequence of this bundle and
Whitehead's theorem  imply the homotopic equivalence of the
manifolds $M$ and $B$. In particular, the fundamental group
$\pi_1(M,x)$ is a free group of rank two.

\vspace{3mm}
{\bf Acknowledgments.} The first author acknowledges the
institutional support of University of Hradec Kr\'alov\'e. The
second author was  supported by the Russian Foundation of Basic
Research (project no. 16-01-00132-a) and by the Program of Basic
Research at the National Research University Higher School of
Economics (project no. 98).

\vskip1cm

A.S. Galaev\\
University of Hradec Kr\'alov\'e, Faculty of Science,\\
Rokitansk\'eho 62, 500~03 Hradec Kr\'alov\'e,  Czech Republic
\\ E-mail:
anton.galaev@uhk.cz

\vskip0.7cm

N.I. Zhukova\\
National Research University Higher School of Economics,
\\ Department of Informatics, Mathematics and Computer
Science,\\ ul. Bolshaja Pecherskaja, 25/12, Nizhny Novgorod,
603155, Russia\\
E-mail: nzhukova@hse.ru

\end{document}